\documentclass{gtart_h}

\def\ifplaintex{\expandafter\ifx\csname documentclass\endcsname\relax}

\def\gtp{{\mathsurround=0pt\it $\cal G\mskip-2mu$eometry \&\ 
$\cal T\!\!$opology $\cal P\!$ublications}}  

\def\recd{{\small Received:\qua\receiveddate\ifx\reviseddate\relax
\else\qquad Revised:\qua\reviseddate\fi\par}} 

\def\lognumber#1{\def\thelognumber{#1}}
\def\volumenumber#1{\def\thevolumenumber{#1}}
\def\volumeyear#1{\def\thevolumeyear{#1}}
\def\papernumber#1{\def\thepapernumber{#1}}
\def\pagenumbers#1#2{\def\startpage{#1}\def\finishpage{#2}}
\def\published#1{\def\publishdate{#1}}

\def\received#1{\def\receiveddate{#1}}

\def\accepted#1{\def\accepteddate{#1}}

\def\asciiauthors#1{\def\theasciiauthors{#1}}
\def\asciiaddress#1{\def\theasciiaddress{#1}}
\def\asciiemail#1{\def\theasciiemail{#1}}

\long\def\asciiabstract#1{\long\def\theasciiabstract{#1}}

\let\\\par\let\thelognumber\relax\let\thevolumenumber\relax
\let\thepapernumber\relax\let\thevolumeyear\relax\let\startpage\relax
\let\finishpage\relax\let\publishdate\relax\let\receiveddate\relax
\let\reviseddate\relax\let\accepteddate\relax\let\theasciititle\relax
\let\theasciiauthors\relax\let\theasciiaddress\relax
\let\theasciiabstract\relax

\let\theasciiemail\relax

\ifplaintex
\font\logobig=cmssbx10 scaled 3836
\font\logomed=cmssbx10 scaled 2557
\else
\font\logobig=cmssbx10 scaled 4200
\font\logomed=cmssbx10 scaled 2800
\fi

\long\def\makeagttitle{   
\count0=\startpage
\agt\hfill      
\hbox to 45truept{\vbox to 0pt{\vglue -13truept{\logomed A\kern -.37em{\logobig 
T}\kern -.38em G}\vss}\hss}
\break
{\small Volume \thevolumenumber\ (\thevolumeyear)
\startpage--\finishpage\nl
Published: \publishdate}

\vglue .25truein

{\parskip=0pt\leftskip 0pt plus
1fil\def\\{\par\smallskip}{\Large\bf\thetitle}\par\medskip} \vglue
0.05truein

{\parskip=0pt\leftskip 0pt plus 1fil\def\\{\par}{\sc\theauthors}
\par\medskip}%
 
\vglue 0.03truein 

{\small\leftskip 25truept\rightskip 25truept{\bf Abstract}\stdspace\theabstract

{\bf AMS Classification}\stdspace\theprimaryclass
\ifx\thesecondaryclass\relax\else; \thesecondaryclass\fi\par
{\bf Keywords}\stdspace \thekeywords\par}\vglue 7truept

}   

\ifplaintex
\font\phead=cmsl9 scaled 950
\font\pnum=cmbx10 scaled 913
\font\pfoot=cmsl9 scaled 950
\headline{\vbox to 0pt{\vskip -4.5mm\line{\small\phead\ifnum
\count0=\startpage ISSN 1472-2739 (on-line) 1472-2747 (printed)
\hfill {\pnum\folio}\else\ifodd\count0\def\\{ }%
\ifx\theshorttitle\relax\thetitle\else\theshorttitle\fi\hfill{\pnum\folio}
\else\def\\{ and }{\pnum\folio}\hfill\ifx\theshortauthors\relax\theauthors
\else\theshortauthors\fi\fi\fi}\vss}}
\footline{\vbox to 0pt{\vglue 0mm\line{\small\pfoot\ifnum\count0=\startpage
\copyright\ \gtp\hfill\else
\agt, Volume \thevolumenumber\ (\thevolumeyear)\hfill\fi}\vss}}
\else
\font=cmsl9 scaled 1050
\font\lnum=cmbx10 
\font=cmsl9 scaled 1050
\makeatletter
\def\@oddhead{{\small\lhead\ifnum\count0=\startpage ISSN 1472-2739 
(on-line) 1472-2747 (printed)\hfill {\lnum\number\count0}\else\ifodd\count0
\def\\{ }\ifx\theshorttitle\relax \thetitle \else\theshorttitle\fi\hfill
{\lnum\number\count0}\else\def\\{ and }{\lnum\number\count0}
\hfill\ifx\theshortauthors\relax 
\theauthors\else\theshortauthors\fi\fi\fi}}\def\@evenhead{\@oddhead}
\def\@oddfoot{\small\lfoot\ifnum\count0=\startpage\copyright\ \gtp\hfill\else
\agt, Volume \thevolumenumber\ (\thevolumeyear)\hfill\fi}
\def\@evenfoot{\@oddfoot}
\makeatother
\fi
\let\maketitlepage\makeagttitle
\let\makeshorttitle\maketitlepage
\let\maketitle\maketitlepage

\newwrite\gtoutfile
\long\gdef\makeheadfile{  
{\def\\{, }\def\s{ }
\immediate\openout\gtoutfile head.xxx
\immediate\write\gtoutfile{Proxy-for: \ifx\theasciiauthors\relax
\theauthors\else\theasciiauthors\fi\s<\ifx\theasciiemail\relax\theemail\else\theasciiemail\fi>}
\immediate\write\gtoutfile{\noexpand\\}
\immediate\write\gtoutfile{Authors: \ifx\theasciiauthors\relax
\theauthors\else\theasciiauthors\fi}
{\def\\{ }\immediate\write\gtoutfile{Title: \ifx\theasciititle\relax
\thetitle\else\theasciititle\fi}}
\immediate\write\gtoutfile{Subj-class: GT or SG, GR etc}
\immediate\write\gtoutfile{MSC-class: \theprimaryclass\ifx\thesecondaryclass\relax\else, \thesecondaryclass\fi}
\immediate\write\gtoutfile{Journal-ref: Algebr. Geom. Topol. \thevolumenumber\s
(\thevolumeyear) \startpage-\finishpage}
\immediate\write\gtoutfile{Comments: Published by Algebraic and
Geometric Topology at}
\immediate\write\gtoutfile{\s\s\s  http://www.maths.warwick.ac.uk/agt/AGTVol\thevolumenumber/agt-\thevolumenumber-\thepapernumber.abs.html}
\immediate\write\gtoutfile{\noexpand\\}
\immediate\write\gtoutfile{}
\ifx\theasciiabstract\relax
\immediate\write\gtoutfile{\theabstract}\else
\immediate\write\gtoutfile{\theasciiabstract}\fi
\immediate\write\gtoutfile{}
\immediate\write\gtoutfile{\noexpand\\}
\immediate\write\gtoutfile{}
\immediate\closeout\gtoutfile}}  

\def\maketitlepage{\makeagttitle\makeheadfile}
\let\makeshorttitle\maketitlepage
\let\maketitle\maketitlepage

\lognumber{36}
\volumenumber{4}
\volumeyear{2004}
\papernumber{36}
\pagenumbers{829}{839}
\received{6 March 2004} 
\accepted{30 September 2004}
\published{7 October 2004}
\input xy
\xyoption{all}

\usepackage{amsmath,amssymb}

\newtheorem{thm}{Theorem}
\newtheorem{lem}[thm]   {Lemma}
\newtheorem{cor}[thm]   {Corollary}

\newtheorem{prop}[thm]  {Proposition}

\newcommand{\term}[1]   {{\it #1}}

\newcommand{\s}     {\Sigma}
\newcommand{\smsh}  {\wedge}
\newcommand{\om}    {\Omega}
\newcommand{\cross} {\times}
\newcommand{\wdg}   {\vee}

\newcommand{\of}    {\circ}

\newcommand{\twdl}  {\widetilde}

\newcommand{\sseq}  {\subseteq}

\newcommand{\Qcat}  {\mathrm{Qcat}}

\renewcommand{\cl}    {\mathrm{cl}}
\newcommand{\cat}   {\mathrm{cat}}

\newcommand{\halfsmash}{\rtimes}

\begin{document}

\title{Implications of the Ganea Condition}
    
\author{Norio Iwase\\Donald Stanley\\Jeffrey Strom}
\shortauthors{Norio Iwase, Donald Stanley and Jeffrey Strom}
\asciiauthors{Norio Iwase, Donald Stanley and Jeffrey Strom}

\address{Faculty of Mathematics, Kyushu University, 
Ropponmatsu 4-2-1\\Fukuoka 810-8560, Japan\\\smallskip\\Department 
of Mathematics and Statistics, University of Regina, 
College West 307.14\\Regina, Saskatchewan, Canada\\\smallskip\\Department 
of Mathematics, Western Michigan University, 
1903 W. Michigan Ave\\Kalamazoo, MI 49008, USA}

\asciiaddress{Faculty of Mathematics, Kyushu University, 
Ropponmatsu 4-2-1\\Fukuoka 810-8560, Japan\\Department 
of Mathematics and Statistics, University of Regina, 
College West 307.14\\Regina, Saskatchewan, Canada\\Department 
of Mathematics, Western Michigan University, 
1903 W. Michigan Ave\\Kalamazoo, MI 49008, USA}

\asciiemail{iwase@math.kyushu-u.ac.jp, 
stanley@math.uregina.ca, jeffrey.strom@wmich.edu}

\gtemail{\mailto{iwase@math.kyushu-u.ac.jp}, 
\mailto{stanley@math.uregina.ca}, \mailto{jeffrey.strom@wmich.edu}}

\begin{abstract}
Suppose the spaces $X$ and $X\cross A$ have the same
Lusternik-Schnirelmann category: $\cat(X\cross A) = \cat(X)$.  Then
there is a {\it strict} inequality $\cat(X\cross (A\halfsmash B)) <
\cat (X) + \cat(A\halfsmash B)$ for every space $B$, provided the
connectivity of $A$ is large enough (depending only on $X$).  This is
applied to give a partial verification of a conjecture of Iwase on the
category of products of spaces with spheres.
\end{abstract}

\asciiabstract{Suppose the spaces X and X cross A have the same
Lusternik-Schnirelmann category: cat(X cross A)= cat(X).  Then there
is a strict inequality cat(X cross (A halfsmash B)) < cat (X) + cat(A
halfsmash B) for every space B, provided the connectivity of A is
large enough (depending only on X).  This is applied to give a partial
verification of a conjecture of Iwase on the category of products of
spaces with spheres.}



\maketitle

\section*{Introduction}\addcontentsline{toc}{section}{Introduction}

The product formula $\cat(X\cross Y) \leq \cat(X) + \cat(Y)$ \cite{Bassi}
is one of the most basic relations of Lusternik-Schnirelmann
category.  Taking $Y= S^r$, it implies
that $\cat(X\cross S^r) \leq \cat(X) + 1$
for any $r > 0$.  In \cite{Ganea1}, Ganea asked whether the inequality
can ever be strict in this special case.   The study of the `Ganea condition'
$\cat(X\cross S^r) = \cat(X) +1$ has been, and remains,  a formidable
challenge to all techniques for the calculation of Lusternik-Schnirelmann
category.  In fact, it was only recently that techniques were developed
which were powerful enough to identify a space which does not satisfy the
Ganea condition \cite{Iwase} (see also \cite{Iwase2,Stanley}).
It is still not well understood exactly which  spaces $X$ do  not 
satisfy the  Ganea  condition,  although it has been 
conjectured that they are precisely those spaces for which 
$\cat(X)$ is not equal to the related invariant $\Qcat(X)$
(see \cite{SST,V}).

Since the failure of the Ganea condition 
appears to be a strange property for a space to have, it is reasonable 
to expect that such failure would have useful and interesting 
implications.   In this paper we explore some of the implications of 
the equation $\cat(X\cross A) = \cat(X)$ for general spaces $A$, and for 
$A= S^r$ in particular.

A brief look at the method of the paper \cite{Iwase}
will help to   put our results into proper perspective.
The  new techniques begin with the following question:
if $Y=X\cup_f e^{t+1}$, the cone on $f:S^t\to X$,
then how can we tell if $\cat(Y)> \cat(X)$?
It is shown (see \cite[Thm.\thinspace 5.2]{Iwase2}
and \cite[Thm.\thinspace 3.6]{Stanley})
that, if $t\geq \dim (X)$,  then 
$\cat(Y)=\cat(X)+1$ if and only if a certain Hopf invariant ${\cal H}_s(f)$
(which is a set of homotopy classes)
does not contain the trivial map $*$. 
It is also shown 
\cite[Thm.\thinspace 3.8]{Iwase2} that if $*\in \Sigma^r{\cal H}_s(f)$, then
$\cat(Y\times S^r)\leq \cat(X)+1$. 
Thus $Y$ does not satisfy Ganea's condition if
$*\not\in{\cal H}_s(f)$, but there is at least one $h\in{\cal H}_s(f)$
such that $\s^r h \simeq *$. 

Of course, if $\s^r h\simeq *$, then 
$\s^{r+1} h\simeq *$ as well, and this suggests the following
conjecture (formulated  in \cite[Conj.\thinspace 1.4]{Iwase}): 

\medskip

\noindent{\bf Conjecture}\qua {\sl  If $\cat(X\cross S^r) = \cat(X)$,
then $\cat(X\cross S^{r+1}) = \cat(X)$.}

\medskip 

\noindent
In this paper we prove that this
conjecture is true, provided $r$ is large enough.  

\begin{thm} \label{thm:conj}
Suppose $X$ is a $(c-1)$-connected space
and let $r > \dim(X) - c\cdot\cat(X) + 2$.
If $\cat(X\cross S^r) = \cat(X)$,
then 
$$
\cat(X\cross S^{t}) = \cat(X)
$$ 
for all $t\geq r$.
\end{thm}

\noindent The conjecture remains open for small values of $r$.

 Our main result is much more general: it shows how the equation
$\cat(X\cross A) =  \cat(X)$ governs the Lusternik-Schnirelmann
category of products of $X$ with a vast collection of other spaces.

\begin{thm}\label{thm:main}
Let  $X$ be a $(c-1)$-connected space and let
$A$ be $(r-1)$-connected with 
$r >  \dim(X) -  c\cdot \cat(X) + 2$.
If $\cat(X\cross A) = \cat(X)$ then
$$
\cat(X\cross (  A  \halfsmash  B )) < \cat(X) + \cat( A \halfsmash B)
$$
for every space $B$.
\end{thm}

Here $A\halfsmash  B = (A\cross B)/ B$ is the 
half-smash product of $A$ with $B$.
When $A$ is a suspension, the half-smash product 
decomposes as $A\halfsmash B \simeq  A \wdg (A\smsh B)$
(see, for example, \cite[Lem.\thinspace 5.9]{Stanley}),
so we obtain the following.

\medskip

\noindent{\bf Corollary}\qua  {\sl
Under the conditions of Theorem \ref{thm:main}, if 
$A$ is a suspension, then 
$$
\cat(X\cross (A\smsh B)) = \cat (X)
$$
for every space $B$.}

\medskip 

\noindent Our partial verification of the conjecture is  
an immediate consequence of this corollary:
it the special case $A= S^r$ and $B = S^{t-r}$.

\medskip

\noindent{\bf Organization of the paper}\qua
In Section 1 we recall the necessary background information 
on homotopy pushouts, cone length and Lusternik-Schnirelmann
category.  We introduce an auxiliary space  
and establish its important properties in Section 2.  
The proof of Theorem \ref{thm:main} is presented in Section 3.

\section{Preliminaries}

In this paper all spaces are based and have the 
pointed homotopy type of CW complexes;
maps and homotopies are also pointed.  We denote by $*$ 
the one point space and any nullhomotopic map.
Much of our exposition uses the language of 
homotopy pushouts; we refer to \cite{Mather} for the 
definitions and basic properties.

\subsection{Homotopy Pushouts}

We begin by recalling some basic facts about homotopy pushout squares.
We call a sequence $A\to B\to C$ a \term{cofiber sequence} if the
associated square
$$
\xymatrix{
A\ar[rr]^f\ar[d] && B\ar[d]
\\
{*}\ar[rr] && C 
}
$$
is a homotopy pushout square.  The space $C$ is called the \term{cofiber}
of the map $f$.  One special case that we use frequently is the 
\term{half-smash product}  $A\rtimes B$, which is the cofiber of the
inclusion $B\to A\cross B$.

%
%
%

Finally, we recall the following result on products 
and homotopy pushouts.

\begin{prop}\label{prop:prod}
Let $X$ be any space.    Consider the squares 
$$
\xymatrix{
A\ar[r]\ar[d] & B\ar[d]\ar@{}[rrrd]|{\txt{and}} 
&&& X\cross A \ar[r]\ar[d] & X\cross B\ar[d]
\\
C\ar[r] &D 
&&&X\cross C\ar[r] &X\cross D.
}
$$
If the first square is a homotopy pushout,
then so is the second.
\end{prop}

\begin{proof}
This follows from Theorem 6.2 in \cite{Steenrod}. 
\end{proof}

\subsection{Cone Length and Category}

A \term{cone decomposition} of a space $Y$ is a diagram 
of the form
$$
\xymatrix{
L_0\ar[d]  & L_1\ar[d] &               & L_{k-1}\ar[d]
\\
Y_0 \ar[r] & Y_1\ar[r] & \cdots \ar[r] & Y_{k-1}\ar[r] & Y_k
\\
} 
$$
in which $Y_0= *$, each sequence $L_i\to Y_i\to Y_{i+1}$ is
a cofiber sequence, and $Y_k \simeq Y$; the displayed 
cone decomposition has \term{length} $k$.  The \term{cone length} 
of $Y$, denoted $\cl(Y)$, is defined by
$$
\cl(Y) = \left\{  
\begin{array}{ll}
0 & \mathrm{if}\ Y\simeq {*}
\\
\infty & \mathrm{if}\ Y\ \mathrm{has\ no\ cone\ decomposition,\  and}
\\
k & \mathrm{if\  the\ shortest\ cone\ decomposition\ of}\  Y\ \mathrm{has\ length}\ k.
\end{array}
\right.
$$

The Lusternik-Schnirelmann category of $X$ may be defined
in terms of the cone length of $X$  by the formula
$$
\cat(X) = \inf\{ \cl(Y)\, |\, X\ \mathrm{is\ a\ homotopy\ retract\ of\ Y} \} .
$$
Berstein and Ganea proved this formula in \cite[Prop.\thinspace 1.7]{B-G}
with $\cl(Y)$ replaced by the strong category of $Y$; the formula above
follows from another result of Ganea ---  strong category is equal to cone length
\cite{Ganea3}.
It follows directly from this definition that if $X$ is a homotopy retract of 
$Y$, then $\cat(X) \leq \cat(Y)$.  The reader may refer to \cite{James}
for a survey of Lusternik-Schnirelmann category.

The category of $X$ can be defined in another way that is 
essential to our work. 
Begin by defining the $0^\mathrm{th}$ Ganea fibration 
sequence 
$\xymatrix@1{F_0(X) \ar[r] &  G_0(X)\ar[r]^-{p_0} & X}$ 
to be the familiar path-loop fibration sequence
$\xymatrix@1{\om (X) \ar[r] &  {\mathcal P}(X)\ar[r]  & X}$. 
Given the $n^{\mathrm{th}}$ Ganea fibration sequence
$$
\xymatrix@1{F_n(X) \ar[r] &  G_n(X)\ar[r]^-{p_n} & X} , 
$$ 
let $\overline{G}_{n+1} (X) = G_n(X) \cup CF_n(X)$ be the cofiber of $p_n$
and define $\overline{p}_{n+1}:\overline{G}_{n+1}(X)\to X$ 
by sending the cone to the
base point of $X$.  The $(n+1)^{\mathrm{st}}$ Ganea fibration 
$p_{n+1}:G_{n+1}(X)\to X$
results from converting the map $\overline {p}_{n+1}$ to a fibration.
The  following result is due to Ganea (cf.  Svarc).

\begin{thm}\label{prop:Ganea}
For any space $X$,
\begin{enumerate}
\item[{\rm (a)}]  $\cl(G_n(X)) \leq n$,
\item[{\rm (b)}]  
the map $p_n: G_n(X)\to X$ has a section if and only if $\cat(X)\leq n$, and 
\item[{\rm (c)}]  
$F_n(X) \simeq (\om(X))^{*(n+1)}$, the  $(n+1)$-fold join of $\om X$ with itself.
\end{enumerate}
\end{thm}

\begin{proof}
Assertion (a) follows immediately from the construction.  
For parts  (b)  and (c), see \cite{Ganea2};
these results also appear, from a different point of view,
in \cite{Svarc}.  
\end{proof}

\section{An Auxilliary Space}

Let   $\twdl G_n$ denote  the homotopy pushout in the square
$$
\xymatrix{
G_{n-1}(X)\ar@{^{(}->}[rr]^-{i_1} \ar[d] && G_{n-1}(X)\cross  A\ar[d]
\\
G_n(X) \ar[rr] && \twdl G_n.
}
$$
The maps $p_n: G_n(X)\to X$ and $1_{A} :  A\to  A$
piece together to give a map $\twdl p_n: \twdl G_n\to X\cross  A$.
The space $\twdl G_n$ and the map $\twdl p_n$
play key roles in the forthcoming constructions;
 this section is devoted to establishing some of their
properties.

\subsection{Category Properties of $\twdl G_n$}

We begin by estimating the category of $\twdl G_n$.

\begin{prop}   \label{prop:upperbound}
For any noncontractible  $A$ and $n>0$, 
$
\cat(\twdl G_n)     <  n  + \cat(A) 
$.
\end{prop}

\begin{proof}
For simplicity in this proof, we write $F_i$ for $F_i(X)$
and $G_i$ for $G_i(X)$.

Let $\cat(A) =k$, so $A$ is a retract of another space $A'$ 
with $\cl(A')= k$.
Let $\twdl G'_n = G_n  \cup G_{n-1} \cross A'$; clearly $\twdl G_n$
is a homotopy retract of $\twdl G_n'$ and so it suffices to show
that $\cl(\twdl G_n')  < n+k$.  Let 
$$
\xymatrix{
L_0\ar[d]  & L_1\ar[d] &               & L_{k-1}\ar[d]
\\
A_0' \ar[r] & A_1'\ar[r] & \cdots \ar[r] & A_{k-1}'\ar[r] & A_k'
\\
} 
$$
be a cone decomposition of $A'$.  We will also use the 
cone decomposition of $G_n$ given by the cofiber sequences
$F_{i-1} \to G_{i-1}\to G_i$.   According to a result of Baues
\cite{Baues} (see also \cite[Prop.\thinspace 2.9]{Stanley2}), 
for each $i$ and $j$ there is a cofiber sequence
$$
\xymatrix@1{ 
F_{i-1} * L_{j-1} \ar[r] & G_i \cross A_{j-1}' \cup G_{i-1} \cross A_j'
\ar[r] & G_i \cross A_j'.
}
$$
Now define subspaces $W_s \sseq \twdl G_n'$
by the formula
$$
W_s = 
\left\{
\begin{array}{ll}
\bigcup_{i+j=s} G_i\cross A_j' 
& 
\mathrm{if}\  s\leq n
\\
 G_{n}\cross A_0'\cup 
\left(
\bigcup_{i+j=s\atop i<n} G_i\cross A_j'
\right)
& 
\mathrm{if}\  s>n 
\end{array}
\right.
$$ 
with the understanding that $A_j' = A_k'$ for all $j\geq k$.
The cofiber sequences guaranteed by Baues' theorem can be 
pieced together with the given cone decompositions of $A'$ and 
$G_n$ to give the cofiber sequences
$$
\xymatrix{
F_s \wdg L_s \wdg
\left(\bigvee_{i+j = s-1\atop i<n-1} F_{i} * L_{j}
\right) \ar[r] &W_s \ar[r] & W_{s+1}
}
$$
for each $s < \min\{ n,k\}$;  when $s\geq n$ we alter the cobase
of the cofiber
sequence by removing the $F_s$ summand,
and when $s\geq k$ we must remove the summand $L_s$.   
Since $\twdl G_n' = W_{n+k-1}$,  we have the result.
\end{proof}

Next, we show  that the   map $\twdl p_n: \twdl G_n\to X\cross   A$ has one of
the category-detecting properties of $p_n: G_n(X\cross  A)\to X\cross   A$.

\begin{prop}\label{prop:getsection}  If  $\cat(X\cross A) = \cat(X) = n$,
then $\twdl p_n$ has a homotopy section.
\end{prop}

\begin{proof}
We follow \cite{CLOT} (see also \cite[Thm.\thinspace 2.7]{Iwase})
and define 
$$
\widehat G_n'(X\cross  A) = \bigcup_{i+j= n} G_i(X) \cross G_j(  A).
$$
There is a natural map $h:\widehat G_n'(X\cross  A)\to X\cross   A$ induced by the 
Ganea fibrations over $X$ and $A$.
According to  \cite[Thm.\thinspace 2.3]{CLOT}, $\cat(X\cross A) = n$
if and only if $h$ has a homotopy section.   

Each map  $G_i(X)\cross G_j(A) \to X\cross A$ (with $j>0$) factors through
$G_i(X) \cross A$ and these factorizations are compatible because $p_{i+1}$
extends $p_i$.  So   
$h$ factors as $\widehat G_n'(X\cross A)\to \twdl G_n\to X\cross A$. 
Therefore, if $\cat(X\cross A) = n$, then  $h$,
and hence $\twdl p_n$, has a section.
\end{proof}

\subsection{Comparison of $\twdl G_n$ with $G_n(X)\cross A$}

Let $j: \twdl G_n\to G_n(X)\cross A$ denote the natural inclusion map.

\begin{prop} 
Assume that $X$ is $(c-1)$-connected and that $A$ is $(r-1)$-connected.
Then the homotopy fiber $F$ of the  map $j$ is  $(nc+r-2)$-connected.
\end{prop}

\begin{proof}
There is a cofiber sequence 
$$
\xymatrix@1{
  \twdl G_n \ar[r]^-j & G_n(X)\cross A\ar[r] &\s F_{n-1}(X) \smsh A . 
}
$$ 
Therefore the homotopy fiber of $j$ has the same connectivity as the space
$\om(\s F_{n-1}(X) \smsh A) \simeq \om(\om(X)^{*n} * A)$, namely $nc + r-2$.
\end{proof}

\begin{cor}\label{cor:inj}
Assume $\dim(Z) < nc + r-2$ and let
$f,g: Z\to \twdl G_n$.   Then $f\simeq g$ if and only if $jf \simeq jg$.
\end{cor}

The proof is standard, and we omit it.

\subsection{New Sections from Old Ones}

Suppose that  $\cat(X) = \cat(X\cross A) = n$.   
By Proposition \ref{prop:getsection}
there is a section $\sigma: X\cross A\to \twdl G_n$ of the   
map $\twdl p_n:\twdl G_n\to X\cross A$.  
Define a new map $\sigma': X  \to G_{n}(X)$ by the diagram
$$
\xymatrix{
X\ar[rrrr]^-{\sigma'} \ar[d]_{i_1} && && G_n(X)  
\\
X\cross A \ar[rr]^-\sigma && \twdl G_n \ar@{^{(}->}[rr]^-j
 && G_n(X) \cross A \ar[u]_{\mathrm{pr}_1}.
}
$$
We need the following basic properties of  $\sigma'$.

\begin{prop} \label{prop:newsection}  
If   $\cat(X\cross A) = \cat(X) = n$, then
\begin{enumerate}
\item[\rm(a)]  $\sigma'$ is a homotopy section of the projection $p_n :G_n(X)\to X$, and 
\item[\rm(b)]
if $X$ is $(c-1)$-connected and  $A$ is $(r-1)$-connected with 
$r>  \dim(X)-nc+2$, then the diagram
$$
\xymatrix{
X\ar[rr]^{\sigma'} \ar[d]_{i_1} && G_n(X) \ar@{^{(}->}[d]^k
\\
X\cross A \ar[rr]^\sigma && \twdl G_n
}$$
commutes up to homotopy.
\end{enumerate}
\end{prop}

\begin{proof} 
First consider the diagram
$$
\xymatrix{
X\ar[rr]^{\sigma'} \ar[d]_{i_1} 
    && G_n(X) \ar@{..>}[d]_k \ar@{=}[rr]
    && G_n(X)  \ar[rr]^{p_n}   
    &&   X
\\
X\cross A \ar[rr]^\sigma \ar[rrrrd]_{1_{X\cross A}}
    && \twdl G_n \ar[rr]^-j
    && G_n(X) \cross A \ar[u]_{\mathrm{pr}_1} 
              \ar[rr]^{\mathrm{pr}_1} \ar[d]^{p_n\cross 1_{A}}
    && G_n(X)\ar[u]_{p_n}\ar[d]^{p_n}  
\\
&&&&
X\cross A\ar[rr]^{\mathrm{pr}_1} && X.
}
$$
The diagram of solid arrows is evidently commutative.
Therefore, we have
$
p_n\of \sigma'  
\simeq \mathrm{pr}_1 \of 1_{X\cross A} \of i_1
\simeq 1_X,
$
proving (a).  

To prove (b) we have to show that two maps $X\to \twdl G_n$
are homotopic.  Since
$
\dim(X) < nc + r -2 
$,
it suffices   by Corollary \ref{cor:inj} 
to show that $j\of (\sigma\of  i_1) \simeq  j\of (k\of \sigma')$. 
Since   $\mathrm{pr}_2\of  j\of (\sigma \of i_1)\simeq *\simeq
\mathrm{pr}_2\of j\of (k\of \sigma')$, it remains
to show that $\mathrm{pr}_1\of j\of (\sigma\of  i_1)
\simeq \mathrm{pr}_1\of j\of (k\of \sigma')$.
But both of these maps are homotopic to $\sigma '$. 
\end{proof}

\section{Proof of the Main Theorem }

\noindent{\bf Proof of Theorem \ref{thm:main}}\qua
We have $n= \cat(X) = \cat(X\cross A)$ by hypothesis.
It follows from Proposition \ref{prop:getsection}
that there is a section $\sigma: X\cross A\to \twdl G_n$ 
of the  map $\twdl p_n:\twdl G_n\to X\cross A$.
We then get the  section $\sigma': X  \to G_{n}(X)$
that was constructed and studied in Section 2.3.
 
Consider the following diagram 
and the induced sequence of maps on the homotopy pushouts of the rows
$$
\xymatrix{
(X\cross A) \cross B \ar[d]_{\sigma\cross 1_B}^{\simeq s}   && 
X\cross B \ar[d]^{\sigma'\cross 1_B} \ar[ll]_{i_1 \cross 1_B}\ar[rr]^{\mathrm{pr}_1}   &&
X \ar[d]^{\sigma'}\ar@{}[rr] &&&
Y\ar[d]
\\
\twdl G_n\cross B \ar[d]_{\twdl p_n \cross 1_B} 
&&
G_n(X) \cross B \ar[rr]^{\mathrm{pr}_1}\ar[ll]_{k \cross 1_B}\ar[d]^{p_n \cross 1_B} 
&&
G_n(X) \ar[d]^{p_n}
& 
{} \ar@{..>}[r]^{\mathrm{homotopy}}_{\mathrm{pushout}}  
&&
P\ar[d]
\\
(X\cross A) \cross B     && 
X\cross B  \ar[ll]_{i_1 \cross 1_B}\ar[rr]^{\mathrm{pr}_1}   &&
X  \ar@{}[rr] &&&
Y.
}
$$
Proposition \ref{prop:newsection} implies that the upper left square commutes
up to homotopy.  
Since $i_1\cross 1_B$ is a cofibration, we can apply homotopy extension
and replace the 
map $\sigma\cross 1_B: (X\cross A)\cross B \to \twdl G_n \cross B$
with a homotopic map $s$ which makes that square strictly commute.
All other squares are strictly commutative as they stand. 

Since the composites 
$(\twdl p_n \cross 1_B) \of (\sigma'\cross 1_B)$
and $p_n\of \sigma'$ are
the identity maps and  $(\twdl p_n \cross 1_B) \of s$
is a homotopy equivalence,   each vertical  composite in the modified diagram is 
a homotopy equivalence.  Thus $Y$ is a homotopy retract of 
$P$,  and consequently $\cat(Y)\leq \cat(P)$.  

The space  $Y$ is   the homotopy pushout of the top row in  the diagram,
which is the   product of the homotopy pushout diagram
$$
\xymatrix{
  B             \ar[rr] \ar[d] &&       {*}           \ar[d]
\\
 A  \cross  B  \ar[rr]    && A\halfsmash B
}
$$
with the space $X$.  Therefore $Y\simeq X\cross (A\halfsmash B)$
by Proposition \ref{prop:prod}.  Since $Y$ is a homotopy retract of 
$P$, it follows that 
$$
\cat(X\cross (A\halfsmash B))  \leq \cat(P),
$$
the proof will be complete once we establish that $\cat(P) < \cat(X) + \cat(A\halfsmash B)$.
This is accomplished in Lemma \ref{lem:P}, which is proved below.
\qed

\begin{lem}\label{lem:P}
The space $P$ constructed in the proof of Theorem
\ref{thm:main} satisfies $\cat(P)\leq \cl(P)
< \cat(X) + \cat( A\rtimes B) $.
\end{lem}

\begin{proof}
The space $\twdl G_n$ is defined by the homotopy pushout square
$$
\xymatrix{
G_{n-1}(X)  \ar[rr]\ar[d]  && G_n(X)  \ar[d]
\\
G_{n-1}(X)\cross A  \ar[rr] && \twdl G_n.
}
$$
Take the product of this square with the space $B$
and adjoin the homotopy pushout square that defines $P$
to obtain the diagram
$$
\xymatrix{   
G_{n-1}(X)\cross B\ar[rr]\ar[d]        && G_n(X)\cross B \ar[rr]\ar[d]   && G_n(X)  \ar[d]
\\
G_{n-1}(X)\cross   A\cross B \ar[rr]  &&\twdl G_n\cross B   \ar[rr]      && P. 
}
$$
By  \cite[Lem.\thinspace 13]{Mather}, the outer  square
$$
\xymatrix{
G_{n-1}(X) \cross   B \ar[rr] \ar[d] && G_n(X) \ar[d]
\\
G_{n-1}(X)  \cross  A \cross   B \ar[rr] && P 
}
$$
is also a homotopy pushout square.  The top map is the composite 
$$
\xymatrix@1{
G_{n-1}(X) \cross B \ar[r]^-{\mathrm{pr}_1} & 
G_{n-1}(X)  \ar@{^{(}->}[r] & G_n(X), 
}
$$ 
and so
we have a new factorization into homotopy pushout squares:
$$
\xymatrix{
G_{n-1}(X) \cross   B \ar[rr]^-{\mathrm{pr}_1} \ar[d] && G_{n-1}(X)  \ar[rr]\ar[d] && G_n(X) \ar[d]
\\
G_{n-1}(X)  \cross A \cross   B \ar[rr]  && L \ar[rr] && P.
}
$$
To identify the space $L$,   observe that the left square is 
simply the product of the space $G_{n-1}(X)$ with the homotopy pushout square
$$
\xymatrix{
  B             \ar[rr] \ar[d] &&      {*}           \ar[d]
\\
 A  \cross  B   \ar[rr]    && A\halfsmash B. 
}
$$
By Proposition \ref{prop:prod},
$L \simeq  G_{n-1}(X)  \cross  (A\halfsmash B)$.
Hence the right-hand square is the  homotopy pushout square
$$
\xymatrix{
G_{n-1}(X)                             \ar[rr]\ar[d]  && G_n(X)    \ar[d] 
\\
G_{n-1}(X)  \cross (A\halfsmash B)  \ar[rr]        && P.
}
$$
Therefore $\cl(P)\leq \cat(X) + \cat(A\halfsmash B)$
by Proposition \ref{prop:upperbound}.
\end{proof}

\Addresses\recd
\end{document}